
\magnification=1100

\def\I#1,#2{I_{#1,#2}}   
\def\II#1,#2{II_{#1,#2}} 
\def\Z{{\bbb Z}}     
\def\Q{{\bbb Q}}     
\def\C{{\bbb C}}     
\def\R{{\bbb R}}		
\def\aut{\mathop{\rm Aut}\nolimits}
\def\mod{\mathop{\rm mod}\nolimits}
\def\H{{\cal H}}	
\def\D{{\cal D}}	
\def\G{\Gamma}		
\def\sset{\subseteq}
\def\Stil{\tilde{S}}	
\def\i{\imath}		
\def\isomorphism{\cong}
\def\w{\omega}		
\def\wbar{\bar{\w}}
\def\tensor{\otimes}
\def\ip#1#2{\langle #1 | #2\rangle}
\def\spanof#1{\langle #1\rangle}
\def\congruent{\equiv}
\def\Ahat{\hat{A}}
\def\Khat{\hat{K}}

\def\rom#1{({\it\romannumeral#1\/})}

\def\qed{\leavevmode\vrule height0pt width0pt
	depth0pt\nobreak\hfill\proofbox\smallskip}
\def\proofbox{\box{1.3ex}{1.3ex}{.1ex}}
\def\box#1#2#3{\hbox{%
\vrule height#2 depth0pt width#3\kern-#3%
\raise#2\hbox{\vrule height0pt depth#3 width#1}\kern-#1%
\vrule height#3 depth0pt width#1\kern-#3%
\vrule height#2 depth0pt width#3%
}}

\newfam\bbbfont
\font\tenbbb=msbm10
\textfont\bbbfont=\tenbbb
\def\bbb{\fam=\bbbfont}
\newfam\calfont
\font\tencal=eusm10
\textfont\calfont=\tencal 
\def\cal{\fam=\calfont}	

\parskip=0pt


\begingroup
\parindent=0pt
{\bf The Period Lattice for Enriques Surfaces}

\medskip

Daniel Allcock\footnote*{Harvard University, Cambridge, MA. 
Supported by an NSF Postdoctoral Fellowship.}

{\it allcock@math.harvard.edu}

14 May 1999

MSC: 14J28 (11F55, 11E12)
\endgroup

\bigskip\noindent
It is well-known that the isomorphism classes of complex
Enriques surfaces are in 1-1 correspondence with a Zariski-open
subset $(\D-H)/\G$ of the quotient of the Hermitian symmetric
space $\D$ for $O(2,10)$.  Here $H$ is a totally geodesic
divisor in $\D$ and $\G$ is a certain arithmetic group. In the
usual formulation of this result [5], $\G$ is described as the
isometry group of a certain integral lattice $N$ of signature
$(2,10)$. This lattice is quite complicated, and requires
sophisticated techniques to work with. The purpose of this note
is to replace $N$ by the much simpler lattice $\I2,{10}$, the
unique odd unimodular lattice of signature $(2,10)$. This allows
for dramatic simplifications in several arguments concerning
$N$, replacing intricate analysis by elementary facts. For
example, in this setting it is easy to see that $H/\G\sset\D/\G$
is irreducible, and also easy to enumerate the boundary
components in the Satake compactification of $\D/\G$. Using
$\I2,{10}$ in place of $N$ also allows one to show that $(\D-H)/\G$
has contractible universal cover. 
The last of these results is new, and the full proof appears in
[1]; here we only give the main idea.
The basis of this paper is a lattice-theoretic trick which 
is well-known to those who work with lattices; however, its
applications in this setting do not appear to have been
published before.

\medskip
We review some notation and facts from [5]. $U$ denotes the
two-dimensional lattice with inner product matrix $\left({0\atop
1}\;{1\atop 0}\right)$. If $A$ is a lattice then $A(n)$ denotes
a copy of $A$ with all inner products multiplied by $n$. An
Enriques surface $S$ has fundamental group $\Z/2$ and its
universal cover $\Stil$ is a K3 surface. The covering
transformation $\i$ acts on $L=H^2(\Stil,\Z)\isomorphism
E_8(-1)^2\oplus U^3$, and its fixed lattice $M$ is primitive in
$L$ and isomorphic to $E_8(-2)\oplus U(2)$. It turns out that
$N=M^\perp$ is isomorphic to $E_8(-2)\oplus U(2)\oplus U$. 
There is a holomorphic 2-form $\w$ on $\Stil$, unique up to a
multiplicative constant, and it satisfies $\i(\w)=-\w$. This
implies that as an element of $H^2(S,\C)=L\tensor\C$, it lies in
$N\tensor\C$. Furthermore, $\w$ satisfies the equality
$\ip{\w}{\w}=0$ and the inequality $\ip{\w}{\wbar}>0$, where
$\ip{\,}{\,}$ denotes the usual bilinear form on $H^2(S,\C)$.

We write $K$ for a fixed copy of $E_8(-2)\oplus U(2)\oplus U$;
by isometrically identifying $N$ with $K$ we may regard $\w$ as
defining an element of $K\tensor\C$. The isometry group of
$K\tensor\R$ is isomorphic to $O(2,10)$, and we write $\D$ for
the Hermitian symmetric space for this group. A concrete model
for $\D$ is the set of points in the complex projective space
$P(K\tensor\C)$ with a representative vector $v$ satisfying
$\ip{v}{v}=0$ and $\ip{v}{\bar v}>0$. That is, under our
identification of $N$ with $K$ we may regard $\w$ as defining an
element of $\D$. This element depends on the choice of
identification, but the point of $\D/\G$ defined by $\w$ does
not, where $\G=\aut K$.

The Torelli theorem for Enriques surfaces asserts that the map
assigning to $S$ the point of $\D/\G$ represented by $\w$ is a
bijection from the the set of isomorphism classes of Enriques
surfaces onto a certain Zariski-open subset $\D_0/\G$ of
$\D/\G$. This map is called the period map. The subset $\D_0$ is
$\D-H$, where
$H=\cup_\ell H_\ell$, $\ell$ varies over the vectors
of $K$ of norm $\ip{\ell}{\ell}=-2$, and $H_\ell$ denotes the set
of points of $\D$ represented by vectors orthogonal to $\ell$.

Our purpose is to restate the Torelli theorem in terms of the
simpler (in particular, unimodular) lattice $\I2,{10}$. For
$p,m\geq0$ we write $\I p,m$ for the lattice with an orthogonal
basis of $p+m$ vectors, $p$ of which have norm $+1$ and $m$ of
which have norm $-1$. If $p,m>0$ then $\I p,m$ is the unique odd
unimodular lattice of signature $(p,m)$. For $p,m>0$ and
$p-m\congruent0\;(\mod 8)$ we write $\II p,m$ for the unique even
unimodular lattice of signature $(p,m)$. This may be described
as a direct sum of various copies of $E_8$, $E_8(-1)$ and
$U$; in particular, $\II1,1\isomorphism U$.
The $\I p,m$ and $\II p,m$ account for all the indefinite
unimodular lattices.

\proclaim Lemma 1.
If  $B$ is even and unimodular and $A\cong B(2)\oplus\II1,1$, then
there is a lattice $\Ahat$ in $A\otimes\R$ isometric to
$B\oplus\I1,1$ such that every isometry of $A$ preserves $\Ahat$
and vice-versa. In particular,
$\aut(B(2)\oplus\II1,1)\cong \aut(B\oplus\I1,1)$.

{\it Proof:} One can check that $(2^{-1/2}A)^*\cong
B\oplus\II1,1(2)$ lies in a unique odd unimodular lattice, which
we take to be $\Ahat$. Here the asterisk denotes the dual
lattice.  By its oddness, unimodularity, and indefiniteness,
$\Ahat$ must be isometric to $B\oplus\I1,1$. One can recover $A$
from $\Ahat$ as $2^{1/2}\cdot(\Ahat^e)^*$ where $\Ahat^e$ is the
even sublattice of $\Ahat$. The intrinsic nature of these
constructions makes it clear that any isometry of $A$ or $\Ahat$
preserves the other.
\qed

\noindent
Applying the lemma to $K=E_8(-2)\oplus U(2)\oplus U\isomorphism
\II1,9(2)\oplus\II1,1$, we find that $\Khat\isomorphism
\II1,9\oplus\I1,1\isomorphism \I2,{10}$. This yields our version
of the Torelli theorem:

\proclaim Theorem 2.
The period map establishes a bijection between the isomorphism
classes of Enriques surfaces and the points of $(\D-H)/\G$,
where $\G=\aut\Khat$ and $H$ is the divisor $\cup_\ell H_\ell$
with $\ell$ varying over the vectors of $\Khat$ of norm $-1$.

{\it Proof:} Applying the lemma we see that $\aut
K=\aut\Khat$. The fact that the description of $H$ given here
coincides with the one given earlier follows from the fact that
the norm $-2$ vectors of $K$ correspond bijectively with the
norm $-1$ vectors of $\Khat$. (Formally, we say that a primitive
vector of $K$ and a primitive vector of $\Khat$ correspond to
each other if they differ by a positive real scalar.) Then our
statement of the Torelli theorem follows from that of Namikawa
([5], thm.~1.14), where all the hard work is actually done.
\qed

Although the lattice $K$ is more closely related to the geometry
of Enriques surfaces than $\Khat$ is, using $\Khat$  
offers some significant advantages. We will illustrate this
with several examples. We will
\rom1
show that $H/\G$ is an irreducible
divisor in $\D/\H$, 
\rom2 
enumerate the boundary components of the
Satake compactification $\overline{\D/\G}$, and 
\rom3
indicate the
proof of the fact that the universal (orbifold) covering space
of $\D_0/\G$ is contractible.

\proclaim Corollary 3.
$H/\G$ is an irreducible divisor of $\D/\G$.

{\it Proof:} This follows from the transitivity of $\aut\Khat$
on the norm $-1$ vectors of $\Khat$. This in turn follows from
the fact that the orthogonal complement in $\Khat$ of such a
vector can only be a copy of $\I2,9$, since it is unimodular of
signature $(2,9)$. Given two norm $-1$ vectors of $\Khat$ it is
now easy to construct an isometry of $\Khat$ carrying one to the
other. 
\qed

\noindent
Corollary~3 was left open by Horikawa [4], and Namikawa ([5],
thm.~2.13) proved it only by relying on a deep theorem of
Nikulin.  Borcherds [2] has found an elementary but still
somewhat involved proof (see the remark after his lemma~2.3).

\proclaim Corollary 4.
There are two orbits of primitive isotropic vectors $v$ in
$\Khat$, one with $v^\perp/\spanof{v}\isomorphism\I1,9$ and one
with $v^\perp/\spanof{v}\isomorphism\II1,9$; these correspond to
the two $0$-dimensional boundary components of
$\overline{\D/\G}$. There are two orbits of $2$-dimensional
primitive isotropic sublattices $V$ in $\Khat$, one with
$V^\perp/V\isomorphism E_8(-1)$ and one with
$V^\perp/V\isomorphism\I0,8$; these correspond to the two
$1$-dimensional boundary components of $\overline{\D/\G}$.

{\it Proof:}
The Satake compactification is defined in terms of the isotropic
sublattices of $\Khat$, so it suffices to prove the transitivity
claims. It is well-known that when $p,m>0$, the orbits of primitive isotropic
vectors in $\I p,m$ are in 1-1
correspondence with the isomorphism classes of unimodular
lattices of signature $(p-1,m-1)$. The lattice corresponding to
a vector $v$ is $v^\perp/\spanof{v}$. The first claim of the
theorem follows because any unimodular lattice of signature
$(1,9)$ is equivalent to either $\I1,9$ or $\II1,9$. We say that
$v$ is odd or even in these cases, respectively.

Now consider a primitive 2-dimensional isotropic sublattice $V$
of $\Khat$. We claim first that $V$ contains an odd vector. For
otherwise it contains an even vector $v$, so that
$v^\perp/\spanof{v}\isomorphism\II1,9$. It is easy to find a
sublattice $A\isomorphism \II1,9$ of $v^\perp$ complementary to
$\spanof{v}$. Since $\Khat$ is odd,
$A^\perp\isomorphism\I1,1$. Taking a primitive vector $w$ of
$A\cap V$ we see that $w^\perp$ contains a copy of $\I1,1$, so
that $w^\perp/\spanof{w}$ is odd, so that $w$ is odd. We have
shown that $V$ contains an odd vector $w$.  Next, consider the
orbits of such pairs $(V,w)$. These are in 1-1 correspondence
with orbits of pairs $(w,W)$ where $w$ is as before and $W$ is a
1-dimensional primitive isotropic lattice in
$w^\perp/\spanof{w}$. Since
$w^\perp/\spanof{w}\isomorphism\I1,9$, there are exactly two
possibilities for $W$ (up to isometry), corresponding to the
unimodular lattices of signature $(0,8)$. In these two cases,
$V^\perp/V\isomorphism\I0,8$ or $E_8(-1)$. Since there are only
two orbits of pairs $(w,W)$, there are two orbits of pairs
$(V,w)$, so there are at most two orbits of sublattices $V$. The
corollary follows.
\qed

\noindent
Corollary~4 was first proven by Sterk in [7], props.~4.5--4.6,
using $K$ instead of $\Khat$. He relied on an intricate analysis
of the orthogonal group of $K$ and also of a certain subgroup of
finite index. He went much further than we have done, by making
a detailed study of the boundary of the Satake compactification
of the quotient of $\D$ by this subgroup.

\proclaim Corollary 5.
The universal cover of $\D_0$ is contractible, as is the
universal orbifold cover of $\D_0/\G$.

{\it Proof sketch:} Since $\D_0$ is an orbifold cover of
$\D_0/\G$, the second claim follows from the first. The first
claim is proven in [1]; the idea is to show that the metric
completion of the universal cover of $\D_0$ is a metric space
with non-positive curvature (in a suitable sense). The essential
ingredient in the proof is that the various components of $H$
meet orthogonally wherever they meet. It is easy to see this in
terms of $\Khat$: if the components $H_r$ and $H_s$ of $H$
meet, where $r$ and $s$ are norm $-1$ vectors of $\Khat$, then $r$
and $s$ span a negative-definite sublattice of $\Khat$. Therefore
$\ip{r}{s}^2<\ip{r}{r}\ip{s}{s}=1$, which requires $\ip{r}{s}=0$,
and we conclude that $H_r$ and $H_s$ meet orthogonally.
\qed

Our final application concerns the norm $-4$ vectors of
$K$. According to [5], thm.~2.15, there are two orbits of such
vectors, distinguished by their orthogonal complements.  Such a
vector $v$ is called even (resp. odd) if $v^\perp$ is isometric
to $E_8(-2)\oplus U\oplus\langle4\rangle$ (resp. to
$E_8(-2)\oplus U(2)\oplus\langle4\rangle$). The vectors of even
type are important because the reflections in them define
isometries of $K$; for example, they play an essential role in
Sterk's work ([7], sec.~4). The odd vectors of norm $-4$ seem to
be less useful. Conveniently, after translating to the $\Khat$
setting the even vectors stand out and the odd vectors are
pushed into the background. That is, the even norm $-4$ vectors
of $K$ correspond bijectively to the norm $-2$ vectors of
$\Khat$, while the odd norm $-4$ vectors of $K$ correspond to
some of the norm $-8$ vectors of $\Khat$.

\medskip
We close with a theorem that is not an application of theorem~2 but rather
a simplification of another part of the arithmetic arguments of
[5]. Part of the background for the Torelli theorem is the
fact that there is essentially only one way to embed
$\II1,9(2)\isomorphism E_8(-2)\oplus U(2)$ as a primitive
sublattice of $\II3,{19}\isomorphism E_8(-1)^2\oplus U^3$. Such
an embedding arises from the map $H^2(S,\Z)\to H^2(\Stil,\Z)$
induced by the covering map from the K3 surface $\Stil$ to the
Enriques surface $S$. The uniqueness theorem implies that the
orthogonal complement of the sublattice is isometric to
$\II1,9(2)\oplus U\isomorphism E_8(-2)\oplus U(2)\oplus U$,
setting the stage for the Torelli theorem. The proof below may
be regarded as a replacement for the argument (2.9) in [5].

\proclaim Theorem 6.
Two primitive sublattices of $\II3,{19}$ that are isometric to
$\II1,9(2)$ are equivalent under the isometry group of
$\II3,{19}$. The orthogonal complement $B$ of such a sublattice
$A$ is isometric to $\II1,9(2)\oplus\II1,1$, and any isometry of
$A$ or $B$ extends to $\II3,{19}$.

{\it Proof:}
First we show that if $A$ is any such sublattice then
$B=A^\perp$ is isometric to $\II1,9(2)\oplus\II1,1$.
Since $\II3,{19}$ is unimodular its images under the
projections to $A\otimes\Q$ and $B\otimes\Q$ are $A^*$ and
$B^*$, and the projections define a bijection (the gluing map) between
$A^*/A$ and $B^*/B$. All elements of $A^*$ have integral norm,
so the same is true of $B^*$. Reducing the norms of elements of
$A^*$ and $B^*$ modulo 2 defines quadratic forms on $A^*/A$ and
$B^*/B$, and since $\II3,{19}$ is even the gluing map is an
isometry. Since $A^*/A$ contains a 5-dimensional isotropic
subspace, so does $B^*/B$, so that $B$ embeds in an even
unimodular lattice $C$, which must of course be isometric to
$\II2,{10}$. Now $C/2C$ is equipped with a nondegenerate
quadratic form obtained by halving norms of elements of $C$ and
then reducing modulo 2. It is obvious that under this form
$B/2C$ is orthogonal to the isotropic subspace
$(2B^*)/2C$. Since $|(2B^*)/2C|=2^5$ and $|B/2C|=2^7$ we see that
$B$ is the preimage in $C$ of the orthogonal complement of a
5-dimensional isotropic subspace of $C/2C$. Since there is only
one such subspace up to isometry of $C/2C$, and every 
isometry of $C/2C$ is induced by one of $C$, $B$ is determined up to
isometry. Since $\II1,9(2)\oplus\II1,1$ is a possibility for $B$ it
is the only possibility. 

Now, the even unimodular lattices containing $A\oplus B$ as a
primitive sublattice are in 1-1 correspondence with the
isometries $A^*/A\to B^*/B$. It is  easy to check
that every isometry of $\II1,9(2)^*/\II1,9(2)$ is induced by one
of $\II1,9(2)$. It follows that every isometry of $A^*/A$
(resp. $B^*/B$) is induced by one of $A$ (resp. $B$). 
Now, if $A'$ is another primitive sublattice of
$\II3,{19}$ isometric to $A$, then its orthogonal complement
$B'$ is isometric to $B$. Choosing any isometry $A\oplus
B\to A'\oplus B'$ carries $\II3,{19}$ to some even unimodular
lattice containing $A'\oplus B'$. By applying a symmetry of $A'$
or $B'$ we may arrange for $\II3,{19}$ to be carried to itself,
so that $A$ and $A'$ are equivalent under $\aut\II3,{19}$.
Finally, the isometries $f$ of $A$ that
extend to $\II3,{19}$ are just those for which there is an
isometry $g$ of $B$ such that the actions of $f$ on
$A^*/A$ and $g$ on $B^*/B$ correspond under the gluing
map. Therefore any symmetry of $A$  extends to $\II3,{19}$. The
same argument also shows that any symmetry of $B$ extends to
$\II3,{19}$. 
\qed

\noindent
Horikawa ([4], thm.~1.5) proved the first part of this theorem by geometric
methods and Namikawa ([5], thm.~1.4) proved it all by using
sophisticated arithmetic results of Nikulin [6]. Our 
argument can be used in more general settings, but our goal in
this paper has been to present the most elementary possible approach
to the arithmetic of the various lattices related to the Torelli
theorem. The theorem is a special case of Nikulin's
theorem~1.14.2, and our proof is presumably a special case of
his proof.

\medskip
{\it Notes:} As mentioned above, the trick used in lemma~1 to
pass between $\II1,1(2)$ and $\I1,1$ is
well-known to experts. 
The passage between $K$ and $\Khat$ is
implicit in the statement of Borcherds ([3], example~13.7) of the
relationship between the automorphic form $\Psi_z$ he constructs
there and the automorphic form $\Phi$ he constructs in [2].  
Also, Sterk ([7], props.~4.3.12 and~4.3.13) explicitly
established a special case of lemma~1, namely the correspondence
between $\I1,9$ and $E_8(-2)\oplus\II1,1$. If he had obtained
this result in its natural generality, it might have simplified
his work.

I am grateful to R.~Borcherds for stimulating discussions
and correspondence.

\bigskip
\noindent
{\bf References}
\medskip

\item{[1]}
D.~J. Allcock.
 Metric curvature of infinite branched covers.
 Preprint 1999.

\item{[2]}
R.~E. Borcherds.
 The moduli space of {E}nriques surfaces and the fake monster {L}ie
  superalgebra.
 {\it Topology}, 35(3):699--710, 1995.

\item{[3]}
R.~E. Borcherds.
 Automorphic forms with singularities on {G}rassmannians.
 {\it Invent. Math.}, 132:491--562, 1998.

\item{[4]}
E.~Horikawa.
 On the periods of {E}nriques surfaces. {I}.
 {\it Math. Ann.}, 234:73--88, 1978.

\item{[5]}
Y.~Namikawa.
 Periods of {E}nriques surfaces.
 {\it Math. Ann.}, 270:201--222, 1985.

\item{[6]}
V.~V. Nikulin. 
 Integral bilinear forms and some of their applications.
 {\it Math. USSR Izv.}, 14:103--167, 1980.

\item{[7]}
H.~Sterk.
 Compactifications of the period space of {E}nriques surfaces part I.
 {\it Math. Z.}, 207:1--36, 1991.

\bigskip
\leftline{Department of Mathematics}
\leftline{Harvard University}
\leftline{One Oxford St}
\leftline{Cambridge, MA 02139}
\smallskip
\leftline{web page: \it http://www.math.harvard.edu/$\sim$allcock}
\bye